# Experimental Study of Nonlinear Resonances and Anti-resonances in a Forced, Ordered Granular Chain

Y. Zhang,[1,a)] D. Pozharskiy,[2] D.M. McFarland,[3]
P.G. Kevrekidis,[4,5] I.G. Kevrekidis,[2,6] and A.F. Vakakis,[1]

[1]*Department of Mechanical Science and Engineering,
University of Illinois at Urbana-Champaign, Urbana, IL 61801, USA*

[2]*Department of Chemical and Biological Engineering,
Princeton University, Princeton, NJ 08544, USA*

[3]*Department of Aerospace Engineering,
University of Illinois at Urbana-Champaign, Urbana, IL 61801, USA*

[4]*Department of Mathematics and Statistics,
University of Massachusetts, Amherst, MA 01003-4515, USA*

[5]*Center for Nonlinear Studies and Theoretical Division,
Los Alamos National Laboratory, Los Alamos, NM 87544, USA*

[6]*Program in Applied and Computational Mathematics,
Princeton University, Princeton, NJ 08544, USA*

## Abstract

We experimentally study a one-dimensional uncompressed granular chain composed of a finite number of identical spherical beads with Hertzian interactions. The chain is harmonically excited by an amplitude- and frequency-dependent boundary drive at its left end and has a fixed boundary at its right end. Such ordered granular media represent an interesting new class of nonlinear acoustic metamaterials, since they exhibit essentially nonlinear acoustics and have been designated as "sonic vacua" due to the fact that their corresponding speed of sound (as defined in classical acoustics) is zero. This paves the way for essentially nonlinear and energy-dependent acoustics with no counterparts in linear theory. We experimentally detect time-periodic, strongly nonlinear resonances whereby the particles (beads) of the granular chain respond at integer multiples of the excitation period, and which correspond to local peaks of the maximum transmitted force at the chain's right, fixed end. In between these resonances we detect a local minimum of the maximum transmitted forces corresponding to an anti-resonance in the stationary-state dynamics. The experimental results of this work confirm previous theoretical predictions, and verify the existence of strongly nonlinear resonance responses in a system with a complete absence of any linear spectrum; as such, the experimentally detected nonlinear resonance spectrum is passively tunable with energy and sensitive to dissipative effects such as internal structural damping in the beads, and friction or plasticity effects. The experimental results are verified by direct numerical simulations and by numerical stability analysis.

---

[a] Corresponding author, yzhng123@illinois.edu

# 1. Introduction

Ordered arrays of granular beads are known to exhibit rich nonlinear dynamical and acoustical behavior, so they have received considerable recent attention in diverse fields of applied mathematics, applied physics and mechanics. This field of research originated from the works of Nesterenko [16] and Lazaridi and Nesterenko [13], and was further explored by other researchers [3, 5, 9, 11, 18, 22-25]. It has been shown that the dynamics and acoustics of uncompressed granular media are highly tunable with energy, and range from being strongly nonlinear and non-smooth in the absence of static pre-compression, to being weakly nonlinear and smooth in the presence of strong pre-compression [5, 15]. Accordingly, depending on the local state of stress, portions of the same granular medium can be either strongly nonlinear (for small stress levels), or weakly nonlinear – even linearized (for stronger states of stress). Two major sources contribute to the nonlinearity of ordered granular media: Nonlinear Hertzian interactions between beads in contact, and bead separations and ensuing collisions in the absence of strong compressive forces. In the limit of no applied pre-compression there is complete absence of linear acoustics in ordered granular media, resulting in zero speed of sound as defined in classical acoustics; accordingly, these media have been characterized as "sonic vacua" by Nesterenko [15]. Despite the absence of linear acoustics, however, it has been shown that these essentially nonlinear sonic vacua possess highly complex nonlinear responses with no counterparts in linear theory, such as nonlinear normal modes [10], moving localized breathers[30], nonlinear traveling waves [26], and frequency bands [6].

Resonance is the main mechanism for energy transmission in spatially periodic systems. Most published works on ordered granular media focus on impulsive inputs which often are simplified as prescribed initial conditions. Only a few works have studied the responses of granular particles to harmonic excitations [2, 4, 8, 14]. Recently, Pozharskiy et al. [20] analyzed the dynamics of a harmonically excited granular medium under a frequency-dependent boundary drive with fixed amplitude. Despite the lack of sound propagation in this medium and the lack of a linear resonance spectrum, various types of purely nonlinear waves were observed involving periodic traveling waves, and a mixture of traveling pulses and standing waves, under different excitation frequencies. Focusing on the force transmitted across the harmonically forced medium, local peaks (resonances) and dips (anti-resonances) of force transmission were noticed, and detailed bifurcation diagrams studying the stability of such responses were reported.

In the present paper we perform an experimental study of a one-dimensional granular medium under harmonic excitation. As this is a strongly nonlinear dynamical system, its dynamics depends critically on the energy level and the system parameters. We numerically show that in this system there occur resonances and anti-resonances in the long-term responses, in the form of traveling and standing waves. We experimentally examine the existence of resonances and anti-resonances by slightly modifying the configuration of the system through the addition of soft flexures. We show that the very rich and complex dynamics of this system can be reproduced by our experimental setup, with good agreement between experimental and numerical results. Furthermore, we examine the stability of the associated solutions by obtaining the corresponding

Floquet multipliers upon converging to the numerically exact (up to a prescribed accuracy) solution via Newton's method. We also quantify the stability of a data-assimilation formulation, based on the experimentally measured force transmitted to the leading bead. In a more general context, this work aims to experimentally prove the realization of strongly nonlinear resonances and anti-resonances in a granular medium with a complete absence of any linear resonance spectrum.

## 2. Strongly nonlinear resonances and anti-resonances in the forced granular medium

A schematic of the considered one-dimensional granular system is presented in Figure 1. This uncompressed granular chain consists of 11 identical particles (beads), and a right fixed boundary. A prescribed harmonic displacement is applied at the left boundary through the (prescribed) motion of a "zero-th" bead. The beads are of spherical shape and are composed of linearly elastic material; moreover, it is assumed that the developed stresses due to bead-to-bead interactions are within the elastic limit of the material of the beads. Under compressive internal forces the interactions between beads obey the essentially nonlinear Hertzian law, with some added dissipative effects as discussed below. In the absence of compressive forces, bead separations and ensuing collisions between them may occur, providing an additional source of strong nonlinearity in the dynamics. To model dissipative effects in the granular chain due to inherent internal structural damping within the beads and frictional effects during bead interactions, we introduce linear damping terms in the model. This type of viscous damping effects has been proven to be adequate in modeling granular dynamics [1, 7, 19, 21]. Moreover, we assume that the beads are constrained to move only in the horizontal direction, so the network is one-dimensional; under this condition additional effects due to rotations of the beads can be neglected, so such rotational effects are not taken into account in the model. Such rotational effects, however, cannot be neglected, e.g., in hexagonally packed granular channels [29].

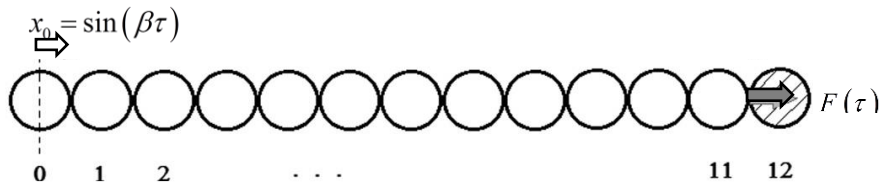

Figure 1. The one-dimensional granular network under periodic excitation.

According to the previous assumptions, the granular medium can be theoretically approximated by a discrete model consisting of concentrated point masses (the beads) with nearest-neighbor coupling stiffnesses obeying the Hertzian (3/2) power force law for compression and exerting zero force (i.e., allowing for separation between beads) in the absence of compression. Accordingly, the equations of motion of the granular system of Figure 1 can be approximately expressed as follows,

$$m\frac{d^2u_1}{dt^2} = \frac{E\sqrt{2R}}{3(1-v^2)}\left\{\left(A_0\sin(2\pi ft)-u_1\right)_+^{3/2} - \left(u_1-u_2\right)_+^{3/2}\right\}$$

$$+D\left\{\left[2\pi fA_0\cos(2\pi ft)-\dot{u}_1\right]\mathrm{H}\left(A_0\sin(2\pi ft)-u_1\right) - \left(\dot{u}_1-\dot{u}_2\right)\mathrm{H}\left(u_1-u_2\right)\right\}$$

$$m\frac{d^2u_i}{dt^2} = \frac{E\sqrt{2R}}{3(1-v^2)}\left\{\left(u_{i-1}-u_i\right)_+^{3/2} - \left(u_i-u_{i+1}\right)_+^{3/2}\right\} \tag{1}$$

$$+D\left\{\left(\dot{u}_{i-1}-\dot{u}_i\right)\mathrm{H}\left(u_{i-1}-u_i\right) - \left(\dot{u}_i-\dot{u}_{i+1}\right)\mathrm{H}\left(u_i-u_{i+1}\right)\right\}, \quad i=2,3,\ldots,N-1$$

$$m\frac{d^2u_N}{dt^2} = \frac{E\sqrt{2R}}{3(1-v^2)}\left\{\left(u_{N-1}-u_N\right)_+^{3/2} - \left(u_N\right)_+^{3/2}\right\}$$

$$+D\left\{\left(\dot{u}_{N-1}-\dot{u}_N\right)\mathrm{H}\left(u_{N-1}-u_N\right) - \left(\dot{u}_N\right)\mathrm{H}\left(u_N\right)\right\}$$

where $m$ is the mass of one spherical bead, $E$ and $v$ the elastic modulus and Poisson's ratio of its material, $R$ its radius (for material density denoted by $\rho$), and $D$ the viscous damping coefficient; the variable $u_i$, $i = 1,2,\ldots,N = 11$, denotes the axial displacement of the $i-th$ bead of the granular chain. The amplitude and frequency of the prescribed displacement excitation are represented by $A_0$ and $f(Hz)$. The subscript $(+)$ denotes that a negative argument in the corresponding parenthesis should be replaced by zero, whereas $H(\cdot)$ denotes the Heaviside function. These terms appear due to the absence of tensile stresses in the chain, implying the possibility of bead separations in the absence of compressive forces. Moreover, zero initial conditions are assumed at the time instant of application of the periodic excitation.

The equations of motion (1) are in dimensional form. For our preliminary computational study reported in this section we consider a granular system composed of steel beads with parameters $m = 28.84 * 10^{-3}\,Kg$, $E = 193 * 10^9\,Pa$, $v = 0.3$, $R = 9.525 \times 10^{-3}m$; and damping parameter $D = 100\,Ns/m$. At this point we would like to emphasize that the system considered in this section is not tied to any specific physical setup since our primary purpose is to highlight the strongly nonlinear resonance and anti-resonance phenomena that can occur in the forced granular chain. Moreover, we assume that the harmonic displacement excitation has constant amplitude $A_0 = 5 \times 10^{-7}m$, which ensures that the elastic deformations of the beads are sufficiently small (within the elastic limit of steel). The dimensional equations of motion can be normalized by dividing each equation by the common mass of the beads, $m$, and introducing the rescalings

$$x_i = u_i/A_0,\ \tau = \varphi t,\ \lambda = D/(m\varphi)\ \text{and}\ \beta = 2\pi f/\varphi$$

where the scaling factor is defined as $\varphi = \left\{\dfrac{\sqrt{2RA_0}\,E}{3(1-v^2)m}\right\}^{1/2}$. Then the non-dimensional equations of motion for this system can be expressed in the following form,

$$\ddot{x}_1 = \left(\sin(\beta\tau) - x_1\right)_+^{3/2} - \left(x_1 - x_2\right)_+^{3/2}$$
$$+ \lambda\left\{\left(\beta\cos(\beta\tau) - \dot{x}_1\right)H\left(\sin(\beta\tau) - x_1\right) - \left(\dot{x}_1 - \dot{x}_2\right)H\left(x_1 - x_2\right)\right\}$$
$$\ddot{x}_i = \left(x_{i-1} - x_i\right)_+^{3/2} - \left(x_i - x_{i+1}\right)_+^{3/2} \qquad (2)$$
$$+ \lambda\left\{\left(\dot{x}_{i-1} - \dot{x}_i\right)H\left(x_{i-1} - x_i\right) - \left(\dot{x}_i - \dot{x}_{i+1}\right)H\left(x_i - x_{i+1}\right)\right\}, \quad i = 2,\ldots,N-1$$
$$\ddot{x}_N = \left(x_{N-1} - x_N\right)_+^{3/2} - \left(x_N\right)_+^{3/2} + \lambda\left\{\left(\dot{x}_{N-1} - \dot{x}_N\right)H\left(x_{N-1} - x_N\right) - \left(\dot{x}_N\right)H\left(x_N\right)\right\}$$

where the variables $x_i$ denote the normalized displacement of the $i-th$ bead of the chain, and $\tau$ is the new normalized time. These equations indicate that the only important parameter governing the nonlinear dynamics is the normalized frequency $\beta$. We study the frequency response of the granular chain by fixing the amplitude of the displacement excitation and recording the maximum of the transmitted force at the rigid wall on the right end for varying frequency. We focus on the stationary-state response of this granular network, i.e., on the state of the dynamics that is eventually reached after a sufficiently long time, so that any initial transients have died out due to dissipative effects. Accordingly, a systematic computational study is performed to investigate the stationary-state nonlinear response of the network under different excitation frequencies.

As an example of the nonlinear stationary-state dynamics of the granular chain [20], in Figure 2 we present the maximum of the transmitted force at the right end of the 11 bead homogeneous chain for varying driving frequency in the range $30Hz \leq f \leq 3000Hz$ and constant amplitude $A_0 = 5*10^{-7} m$. In this particular example there are five distinguishable peaks or "resonances," corresponding to frequencies where maximal force is transmitted to the right end of the chain. The diagram of Figure 2 also features dips or "anti-resonances" between the resonances, associated with frequencies of minimal force transmission. In recent studies it was theoretically [10] and experimentally [6] shown that ordered granular media such as the granular chain of Figure 1 possess energy-tunable pass and stop bands similar to linear spatially periodic systems. These strongly nonlinear frequency bands are due to the essentially nonlinear character of the sonic vacuum and correspond to frequency and energy ranges where disturbances will either propagate (pass bands) or attenuate (stop bands) in the medium. Indeed, in low-frequency acoustic pass-bands, these highly discontinuous media support solitary-like pulses or spatially extended wave transmission. Referring to the plot of Figure 2, the frequency range (coinciding with the pass band) where these resonances and anti-resonances occur lies within the nonlinear pass band of the granular chain of Figure 1, whereas for frequencies above 2000Hz (in the stop band) there is negligible force transmission due to the incapacity of the granular crystal to transmit energy at high frequency. Focusing on the responses within the pass bands where resonances and anti-resonances occur, we explore these features in detail to fully understand the strongly nonlinear mechanisms governing these interesting stationary-state phenomena.

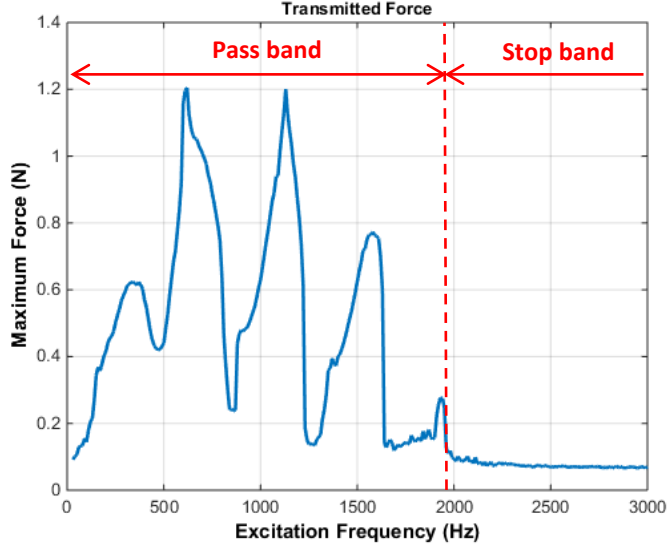

Figure 2. Maximum force transmitted on the right fixed boundary (Pozharskiy et al., 2015).

As a representative example of the dynamic response of the forced granular chain in a resonance, in Figure 3 we depict the spatiotemporal variation of the kinetic energies of the 11 beads at the first (cf. Fig. 3a) and third (cf. Fig. 3b) resonance peaks. It is apparent that the response of the chain is in the form of waves that are initiated on the left (forced) end of the chain, and reflected at the right fixed end. Focusing at the resonance depicted in Fig. 3a we note that the zero-th bead exerts a strong impulse-like excitation at precisely the time instant when the traveling pulse reflected from the right, immovable boundary reaches the left end; one deduces that there is a strong excitation exerted at the granular chain at each period of the prescribed displacement excitation, so that the traveling pulses in the chain fully synchronize with the excitation source. Accordingly, this peak is designated as a 1:1 resonance [20].

On the contrary, at the resonance depicted in Fig. 3b (corresponding to the third peak of Fig. 2), one notes that there is one strong impulse-like excitation every three periods of the zero-th bead oscillation, since the time needed by the traveling pulse to fully traverse twice the length of the chain is three times the period of the excitation; hence, this is designated as a 1:3 resonance. We note that whereas secondary, weaker impulses also occur in between the strong impulse excitations in this case, these are too weak to initiate new identifiable traveling pulses in the granules. Using this classification the resonance peaks of Fig. 2 can be classified as $1:n$ resonances for $n = 1, \ldots, 5$, whereas higher-order resonances are eliminated due to dissipative effects. Moreover, all these nonlinear resonances occur in the pass band of the granular medium of Fig. 1, since only in that frequency range can travelling pulses propagate from the left to the right boundary and vice versa.

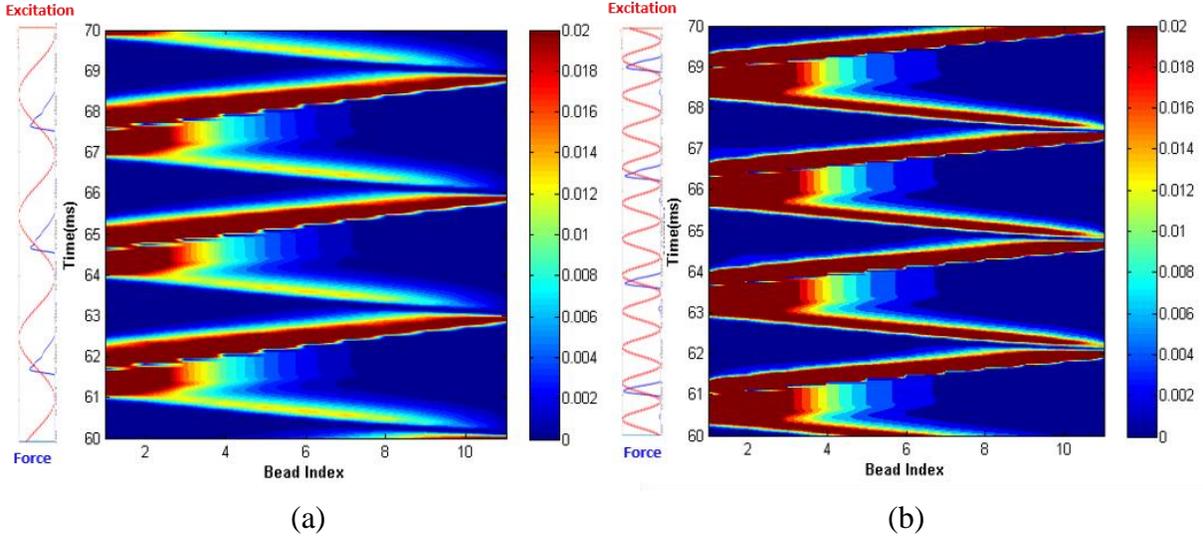

Figure 3. Stationary spatiotemporal variation of the kinetic energy (scale $\times 10^6$) for (a) 1:1 resonance (340 Hz) and (b) 1:3 resonance (1130 Hz); the input harmonic displacement (red) and applied force excitations (blue) for each case are also shown on the left of each plot.

In anti-resonances there occurs strong attenuation of the propagating pulses caused by destructive interference between right- and left-going pulses propagating through the granular chain [20]. Moreover, all anti-resonances are located within the pass band, and correspond to 1:1 synchronization between the prescribed amplitude excitation and the resulting pulses in the medium. Additional unstable or quasi-periodic responses were identified in the theoretical work by Pozharskiy et al. [20], but these are not pursued further herein. In the next section we describe the results of an experimental study that was performed in order to experimentally prove the existence of the theoretically predicted resonances and anti-resonances in a practical granular chain with the configuration shown in Figure 1.

## 3. Experimental study

The experimental study aims to verify the previous theoretical results by means of a harmonically forced granular fixture. For reasons of practical implementation we focus only on low-frequency regimes where dissipative effects are less pronounced, and where low-order resonances and anti-resonances are realized. The experimental fixture is shown in Figure 4, and represents the practical realization of the forced granular medium depicted in Figure 1. It consists of two sturdy pillars connected through threaded shafts. To host the flexures, stainless steel holders with slots are placed on these shafts. Each bead in the granular chain is rigidly attached to the one end of a thin steel flexure and aligned horizontally and vertically. This alignment is crucial in order to experimentally realize the one-dimensionality of the granular dynamics and minimize frictional effects due to relative rotations between adjacent beads. The other end of the aforementioned flexure is placed in a slot of the holders assembled on the threaded shafts. The thin flexures, made of spring steel grade 1095, are designed to be much softer than the stiff beads, so that the time scales of the dynamics of the flexure responses and the dynamics of the bead-to-bead interactions (through Hertzian contact) are separable. It follows that the dynamics of the flexures should minimally affect the granular dynamics, so their contributions to the measured responses should be small. Nevertheless, as discussed below, the theoretical model needs to be modified by weak grounded stiffness in order to account for "ringing effects" in the measured responses due to the dynamics of the supporting flexures. After completing the alignment of the granular chain the holders are firmly bolted and rigidly fastened to the support structures.

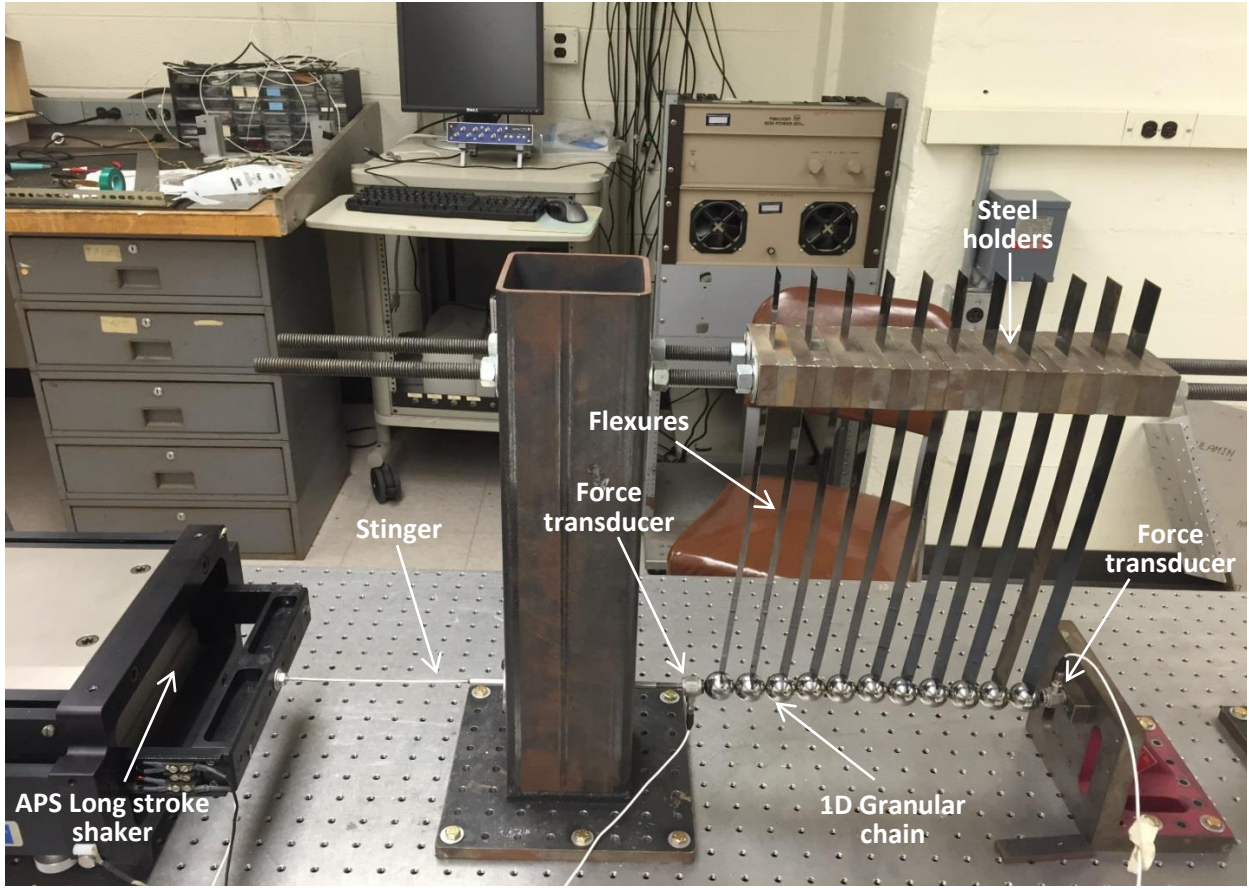

Figure 4. Experimental fixture of the forced granular chain.

The homogeneous granular chain consists of 11 spherical granules composed of bearing-quality aircraft-grade E52100 alloy steel of common radius $R = 12.7\ mm$, modulus of elasticity $E = 210\ GPa$, density $\rho = 7850\ Kg/m^3$ and Poisson's ratio $\nu = 0.3$. The supporting flexure is inserted to a depth of about $1/8^{th}$ of its diameter and permanently glued. Controlled excitation to the first particle (at the left end) is provided by means of an APS® long-stroke shaker. The stinger of the shaker is guided to excite the chain horizontally, and a piezoelectric force transducer (PCB® model 208C03, with sensitivity 2,248mV/kN) is attached at the point of contact of the stinger with the sample in order to measure the applied force. Due to the strong nonlinear bead-stinger dynamic interaction during the measurement, the measured applied force is affected by the measured response since the force sensor is not glued to the first particle, which raises the possibility of losing contact with the chain as we will see below. Hence, a laser vibrometer (Polytec® model PSV-300-U) is employed here to record the velocity of the armature of the shaker which is unaffected by the measured dynamics and enables accurate measurement of the amplitude and frequency of the applied harmonic motion. As a measured output signal, the transmitted force at the right end of the chain is recorded by an additional piezoelectric force transducer (PCB® model 208C02, with sensitivity 11,241mV/kN). To summarize, three measurements (applied force, input velocity and output force) are recorded in the testing of this granular chain. The data is then post processed using Matlab®.

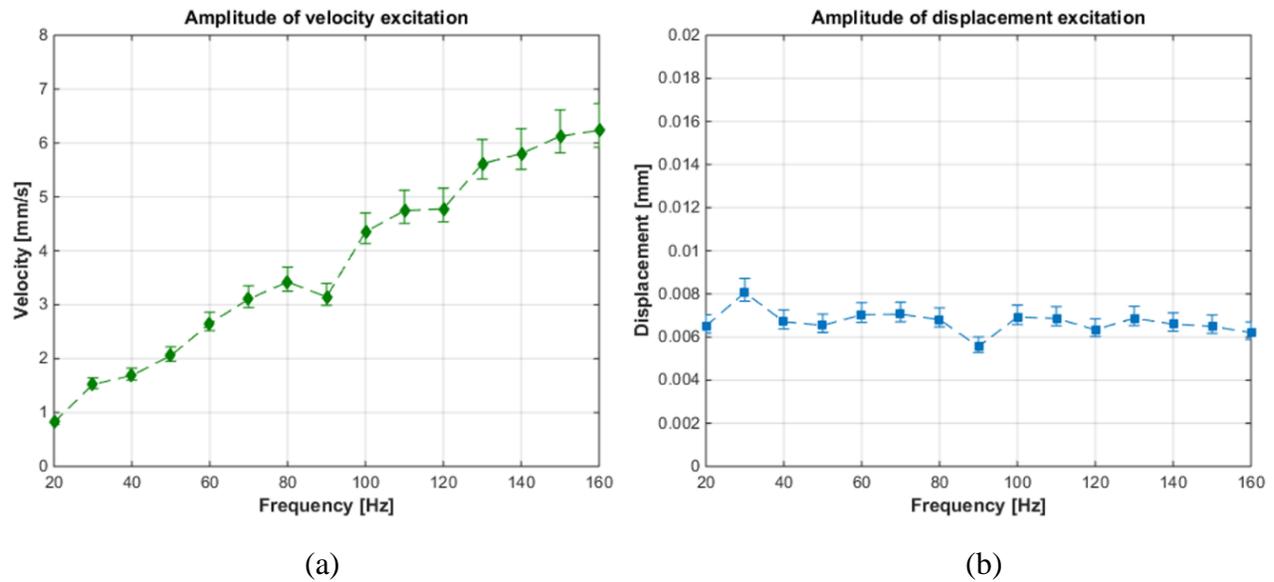

(a)                                                    (b)

Figure 5. Experimental shaker output over the frequency range of the tests: (a) Measured velocity amplitude, and (b) derived displacement amplitude of the armature of the shaker.

As mentioned above, our experimental study considers the forced dynamics of the granular chain under low-frequency harmonic excitation. We aim to maintain the shaker's displacement amplitude as constant as the frequency varies, but this amplitude cannot be recorded directly due to the limitation of the experimental setup. Assuming that the motion of the shaker is ideally harmonic and its initial phase is zero, its velocity can be expressed as $V(t) = V_0\cos(2\pi f t)$, and the corresponding displacement as $A(t) = A_0\sin(2\pi f t)$ where $f$ is the driving frequency and $A_0 = V_0/2\pi f$. In Figure 5 we depict the experimental measurements of the amplitude of the velocity of the shaker ($V_0$) in the range $20Hz \leq f \leq 160Hz$, along with the corresponding displacement amplitude of the shaker ($A_0$). Clearly, the displacement amplitude (cf. Figure 5b) remains nearly constant under different drive frequencies. By maintaining the amplitude of the excitation to a nearly constant level, the only factor that affects the response of tested system is the driving frequency. Hence, all experimental measurements can be plotted and analyzed as functions of the external frequency.

The study of the force transmitted by the 11th (farthest to the right) bead to the force transducer at the right boundary allows us to identify the resonances and anti-resonances in the stationary-state dynamics: Local peaks of the maximum transmitted force correspond to maxima of energy transmission through the granular medium, whereas, local valleys of the maximum transmitted force indicate weak energy transmission. This method for detecting the nonlinear resonances in the forced granular medium was theoretically confirmed by the computational analysis of Pozharskiy et al. [20]. The central result of the experimental study is summarized in Figure 6a, depicting the maxima of the experimentally measured transmitted force at the right of the granular chain for varying drive frequency in the range $20 Hz \leq f \leq 160 Hz$. Within this range we are able to detect two clear local peaks or resonances (located at $60 Hz$ and $100 Hz$), and one valley or anti-resonance in between the two peaks (located at $90 Hz$). Based on the theoretical study of Pozharskiy et al. [20] we identify that these two peaks of transmitted force are 1:1 and 1:2 resonances, respectively, whereas the valley in between is a 1:1 anti-resonance. To ensure repeatability of the experimental results, three different tests were repeated at every frequency, except for the two resonance peaks and the anti-resonance valley, for which five repeated tests were performed at each frequency. Even though the amplitude of the excitation was maintained almost constant in each of the repeated tests, some small variation was still unavoidable during the experiments. Accordingly, the averages of the different trial tests were computed to determine the final measured responses for varying frequency. The extreme measured values of these tests were also indicated in Figures 5 and 6 through the indicated error bars.

A special note is warranted at this point concerning the force excitation exerted by the stinger of the shaker to the first bead of the granular chain. As discussed in previous works [6, 20], at the low-frequency pass band, even though the stinger has a prescribed harmonic motion, it does not maintain continuous contact with the first bead of the granular chain; as a result, the force excitation applied to the granular chain consists of a periodic or quasi-periodic series of force pulses as shown in the time series of input force in Figure 6b. In the pass band each applied force pulse generates a transmitted pulse in the granular chain (cf. Fig. 3a,b) which propagates almost unattenuated (except for dissipative effects) through the chain. This is verified by the force pulses transmitted to the right end of the chain and recorded by the force transducer at the right end (cf. Fig. 6b).

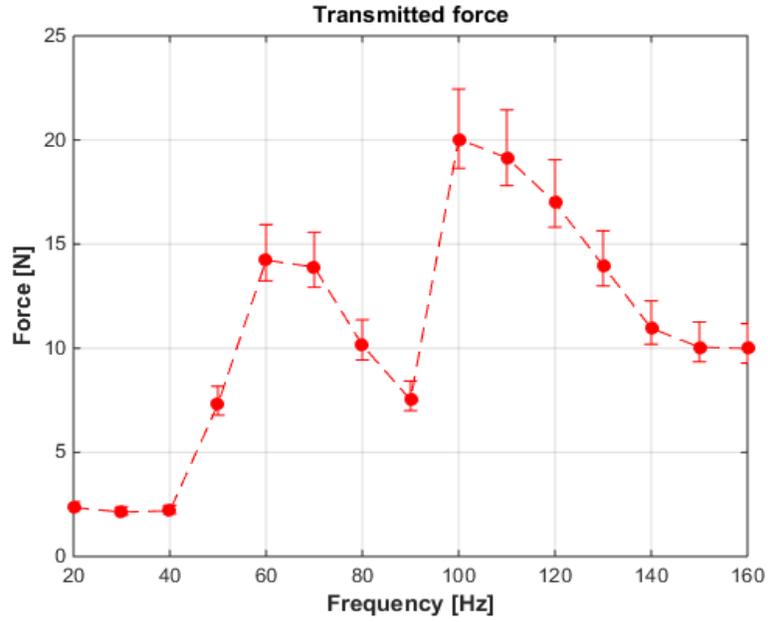

(a)

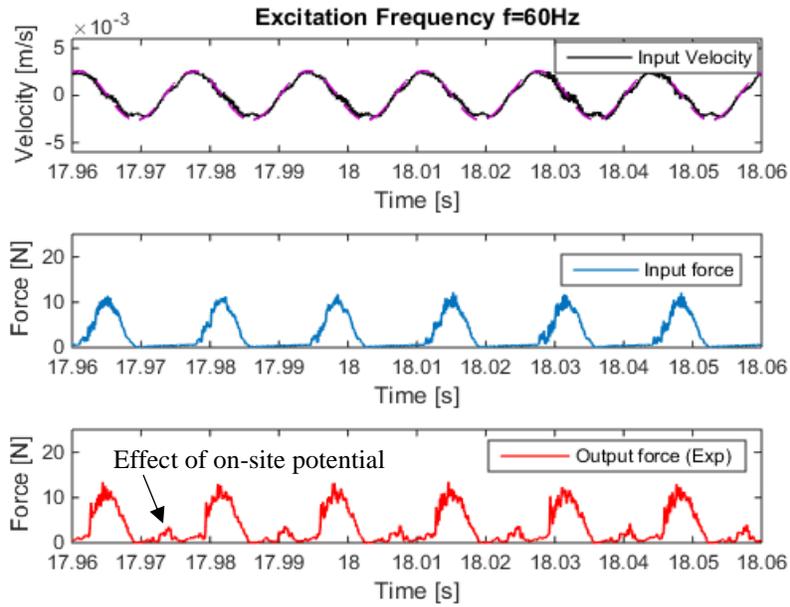

(b)

Figure 6. Experimental response of the granular chain under harmonic excitation: (a) Maximum of the transmitted force measured by the force transducer at the right end for $20Hz \leq f \leq 160Hz$; (b) experimental velocity time series of the armature of the shaker measured by laser vibrometry (*top*), time series of input force applied at the first bead measured by the force transducer at the left end (*middle*), and time series of the transmitted force by the 11$^{th}$ bead measured by the force transducer at the right end (*bottom*), under harmonic excitation at 60 Hz (1:1 resonance).

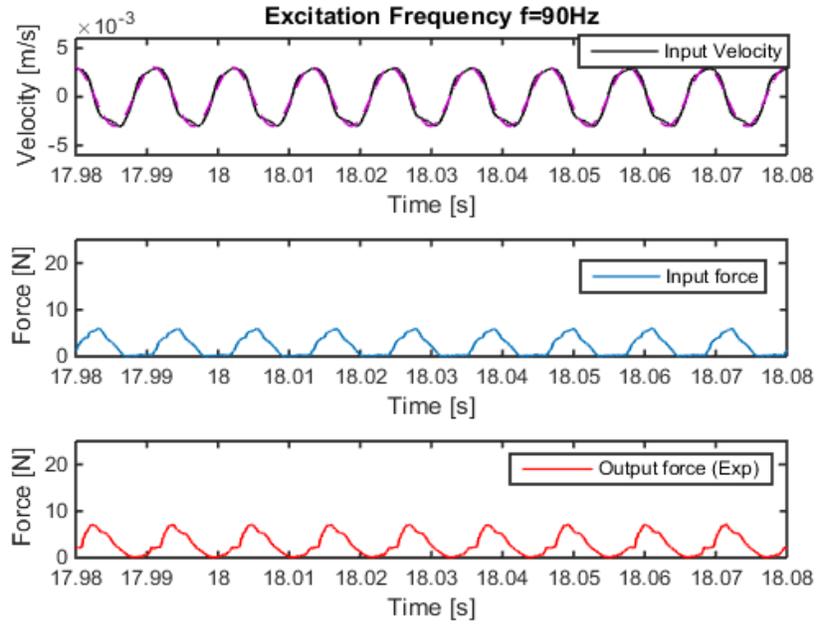

(a)

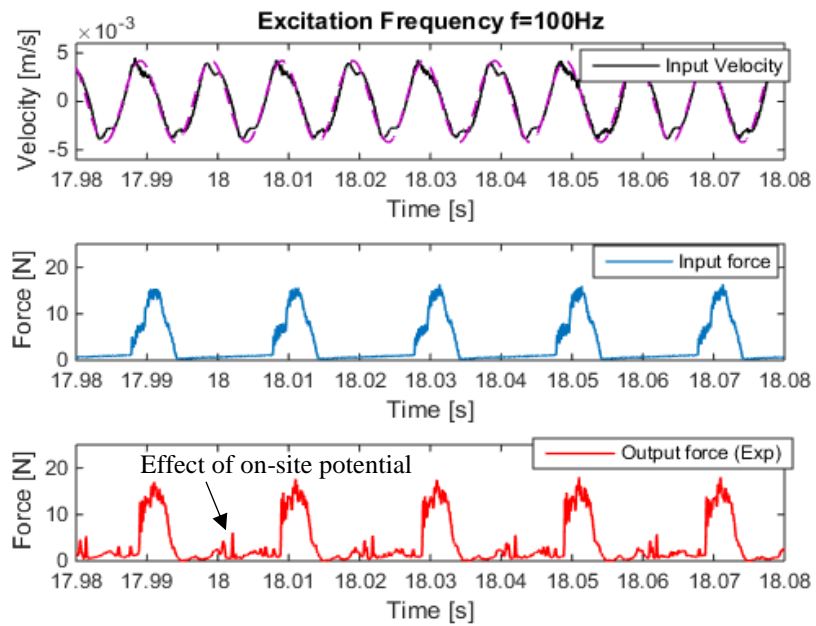

(b)

Figure 7. Experimental velocity time series of the armature of the shaker measured by laser vibrometry (*top*), time series of input force applied at the first bead measured by the force transducer at the left end (*middle*), and time series of the transmitted force by the 11$^{th}$ bead measured by the force transducer at the right end (*bottom*), under harmonic excitation at (a) 90 Hz (1:1 anti-resonance), and (b) 100 Hz (1:2 resonance).

In Figures 6b and 7a,b we depict experimentally measured time series for three distinct forcing frequencies, namely $60 Hz$, $90 Hz$ and $100 Hz$, corresponding to 1:1 resonance, 1:1 anti-resonance, and 1:2 resonance, respectively. In particular, we depict the velocity time series of the armature of the shaker measured by a non-contacting laser vibrometer with reference to numerical harmonic signals at corresponding frequencies; this time series is nearly harmonic, verifying that the applied displacement excitation of the shaker stinger is harmonic. In addition, we depict the time series of the input force applied by the stinger to the first bead of the granular chain and measured by the force transducer located at the left end of the chain; the time series verifies the pulse-like excitation of the granular chain in the low-frequency pass band, as discussed previously. Finally, the time series of the force transmitted by the $11^{th}$ (furthest to the right) bead to the force transducer at the right end of the chain; we note that the transmitted force is again in the form of pulses since each applied impulse pulse on the left end gives rise to a transmitted force pulse at the right end. The recording of the maxima of the transmitted force pulses for varying frequency generates the transmitted force plot of Figure 6a. Clear "silent zones" between any two successive input force pulses are due to the separation between the stinger and the first bead.

However, a different result is observed considering the region between two successive transmitted force pulses in the two resonance peaks (cf. Figs. 6b and 7b) since instead of silent zones we note the realization of distinct small peaks of transmitted force between the main transmitted pulses. These small residual force pulses (or "ringing") are mainly contributed by the presence of the on-site potential (elastic foundation) in the nonlinear dynamics due to the fact that each of the beads of the granular chain is supported (grounded) by the "soft" steel flexures. Even though there is time scale separation between the "soft" dynamics of the flexures and the "stiff" dynamics of the bead-bead interactions, there is still a measured dynamic effect contributed by the supporting flexures. As a result, the $11^{th}$ bead does not fully relax after each interaction with the force transducer (which generates the main transmitted force pulse), but due to the restoring effect of the flexure foundation the bead regains contact with the force transducer, resulting in the small force pulse which occurs after the transmission of the main force pulse. We note, however, that this ringing residual effect is small compared to the main transmitted force pulse train, a result that verifies the adequate time scale separation of the stiff/soft dynamics.

In our study we focus only on the main transmitted forced pulses caused by the Hertzian interactions between beads. Considering first the case of excitation frequency at 60 Hz, there occurs 1:1 resonance in the dynamics (cf. Figs. 6a,b). At this resonance peak the stinger of the shaker exerts a strong force pulse at every period of the harmonic oscillation of the shaker armature, so there is 1:1 synchronization between the displacement of the shaker and the resulting input force to the granular chain. As explained theoretically in Pozharskiy et al. [20], the cause of this 1:1 resonance lies in the relative phase between the periodic motions of the stinger and the first bead of the chain, which is equal to π/2 (condition for resonance) at the instant of contact, at every cycle of the shaker motion. Hence, a strong force pulse is exerted on the chain by the stinger of the shaker at every cycle of the prescribed motion of the shaker, which, in turn gives rise to a similarly synchronized transmitted force pulse at the right end of the chain.

By increasing the excitation frequency to 90 Hz (1:1 anti-resonance), a different type of forced dynamics of the granular chain is identified, where the maximum transmitted force plot reaches a local valley (cf. Figs. 6a and 7a). At that frequency input force pulses delivered by the excitation source have a 1:1 correspondence with low-intensity transmitted force pulses, due to destructive interference of left- and right-going pulses in the granular chain.

Interestingly enough, by further increasing the frequency to 100 Hz we notice a second resonance peak in the transmitted force plot which corresponds to an 1:2 resonance; in this case the granular chain undergoes a period-2 subharmonic motion. Indeed, at this resonance peak (cf. Figs. 6a and 7b) the excitation source exerts a strong force pulse at every second period of the harmonic oscillation of the shaker armature. As a result, the input force pulse measured at the left end of the granular chain repeats itself every two periods of the prescribed harmonic oscillation of the armature of the shaker. Between two successive input force pulses clear silent zones are noticed due to the loss of contact between the stinger and the first bead of the granular chain. As explained theoretically in Pozharskiy et al. [20], the cause of this 1:2 resonance (similarly for all higher order resonances) is the relative phase between the periodic motions of the stinger and the first bead of the chain which leads to loss of contact at every period of the prescribed displacement excitation of the shaker armature; rather, following one cycle of oscillation (period) after a strong force pulse is exerted on the chain, the stinger motion becomes out-of-phase with respect to the motion of the first bead of the chain, so no contact is possible at that time instant. However, after another cycle, the motion of the shaker and of the first bead are in phase difference of $\pi/2$ (condition for resonance) so a strong force impulse is exerted again by the stinger of the shaker to the chain. Moreover, similar to the case of 1:1 resonance (cf. Fig. 6b), small residual force pulses are introduced in the transmitted force measurement due to the soft dynamics of the flexures (cf. Fig. 7b).

Hence, our experimental results fully verify previous theoretical predictions and confirm the existence of strongly nonlinear resonance motions in the granular chain of Figure 4. We note at this point that since the granular chain has no prior compression, it represents a sonic vacuum with complete lack of any linear resonance spectrum. It follows that the experimentally measured resonance spectrum is strongly nonlinear and fully tunable with energy, so the detected resonances and anti-resonances are highly sensitive to the applied energy input. Moreover, due to the highly complex granular dynamics the nonlinear resonance spectrum is highly sensitive to damping, since for low-enough dissipation the granular chain possesses chaotic dynamics (due to separations and ensuing collisions between beads) so no detectable resonance spectrum can be observed. This was the case of a previous study [14] where chaotic resonance motions were experimentally detected in a two-bead granular system, and no discernable resonance spectrum similar to the one presented here (cf. Fig. 6a) could be experimentally detected. This was due to sensitive dependence on initial conditions of the chaotic motions of that system. Chaotic dynamics was also observed systematically in distributed versions of relevant driven-damped granular chains e.g., in [2, 8].

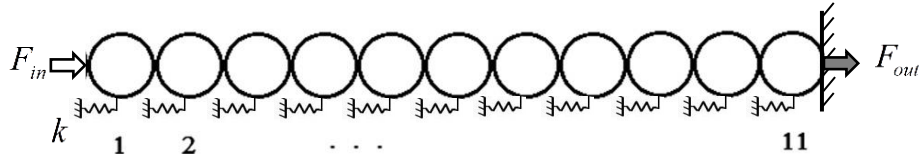

Figure 8. Augmented mathematical model for the forced granular chain incorporating the weak on-site foundation generated by the supporting flexures.

As a further step we aim to reconstruct computationally the experimental measurements. Based on the previous discussion it is clear that the mathematical model (1-2) should be augmented by stiffness terms in order to model the "soft" dynamics of the supporting flexures that give rise to the "ringing" effects in the transmitted force time series of Figs. 6b and 7b. To this end, the augmented configuration depicted in Figure 8 is considered, which compared to the model of Figure 1, possesses weak grounding linear stiffnesses for each bead of the granular chain. Maintaining the notation of the model (1-2), the augmented equations of motion (in dimensional form) of this system are expressed as:

$$m\frac{d^2 u_1}{dt^2} = F_{in} - \frac{E\sqrt{2R}}{3(1-v^2)}\left\{(u_1 - u_2)_+^{3/2}\right\} + D\left\{-(\dot{u}_1 - \dot{u}_2)H(u_1 - u_2)\right\} - k u_1$$

$$m\frac{d^2 u_i}{dt^2} = \frac{E\sqrt{2R}}{3(1-v^2)}\left\{(u_{i-1} - u_i)_+^{3/2} - (u_i - u_{i+1})_+^{3/2}\right\}$$
$$+ D\left\{(\dot{u}_{i-1} - \dot{u}_i)H(u_{i-1} - u_i) - (\dot{u}_i - \dot{u}_{i+1})H(u_i - u_{i+1})\right\} - k u_i, \quad i = 2, 3, \ldots 10 \quad (3)$$

$$m\frac{d^2 u_{11}}{dt^2} = \frac{E\sqrt{2R}}{3(1-v^2)}\left\{(u_{10} - u_{11})_+^{3/2} - \sqrt{2}(u_{11})_+^{3/2}\right\}$$
$$+ D\left\{(\dot{u}_{10} - \dot{u}_{11})H(u_{10} - u_{11}) - (\dot{u}_{11})H(u_{11})\right\} - k u_{11}$$

The first term, $F_{in}$, on the right-hand-side of the equation of motion of the first bead is taken directly from the experimental measurements (cf. Figs. 6b, 7a and 7b), in order to ensure that the computational model is excited by the same pulse train excitation exerted in the experimental fixture. This turns our modeling into a data assimilation problem, since the response of the system is driven by an experimentally obtained signal; the model predictions will then become nonlinear observations of the experimental system state. The numerical values of all parameters in the model (3) were previously defined, except for the damping coefficient $D$ and the (soft) foundation stiffness $k$. The viscous damping in a one dimensional homogenous chain composed of steel beads has been estimated by Herbold et al. [7] as 32.15 Ns/m, and by Potekin et al. [19], as 35.4 Ns/m. Regarding the numerical value of the damping coefficient, we consider the estimated value by Potekin et al. [19] which was derived for the identical experimental setup. Accordingly, we set D = 35.4 Ns/m and use this numerical damping value in all of the following computational simulations that are performed for comparing with the experimental results. To simulate the weak foundation effect resulting from the thin steel flexures, we choose $k$ as 0.1% of the Hertzian contact stiffness between particles, and set $k = E\sqrt{2R}/3000(1-v^2)$.

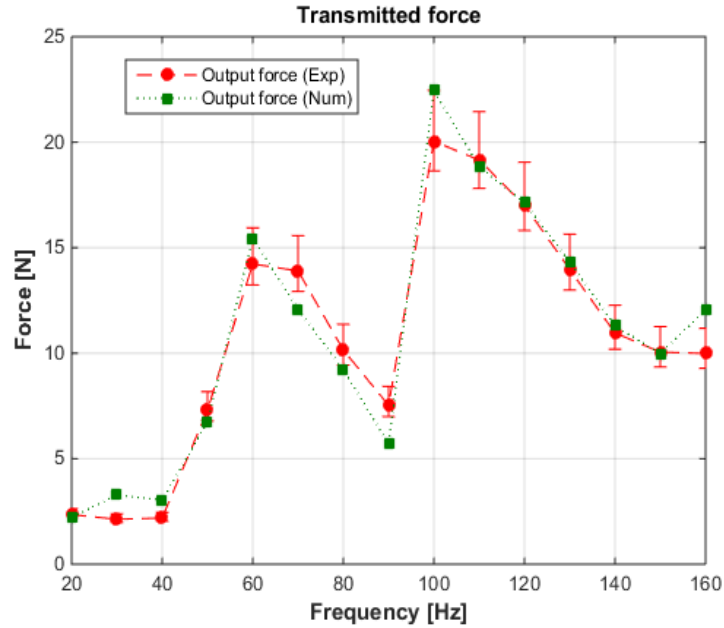

(a)

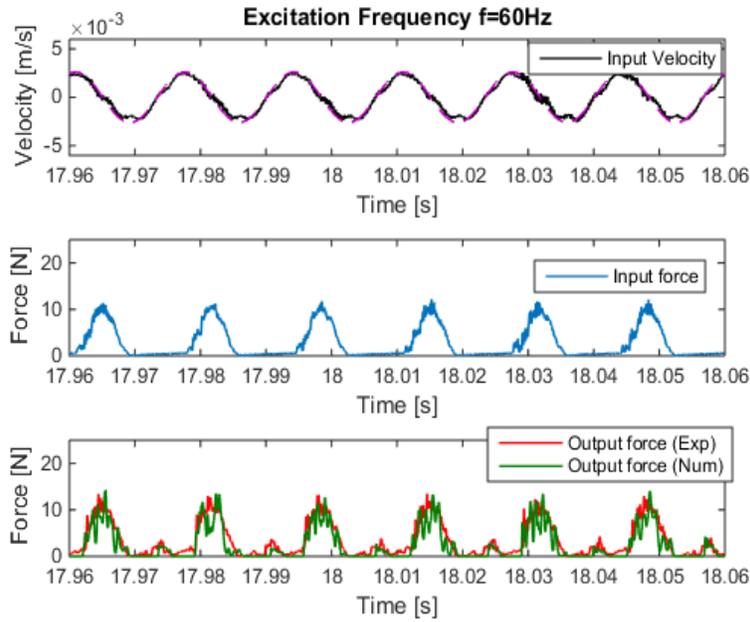

(b)

Figure 9. Comparisons between experimental measurements (red) and numerical results (green) from the model (3): (a) Maximum of the transmitted force at the right end of the granular chain in the frequency range $20 Hz \leq f \leq 160 Hz$; (b) experimental velocity time series of the armature of the shaker (*top*), experimental input force at the left end of the first granule (*middle*) and transmitted force at the right end (*bottom*) for 1:1 resonance at 60 Hz.

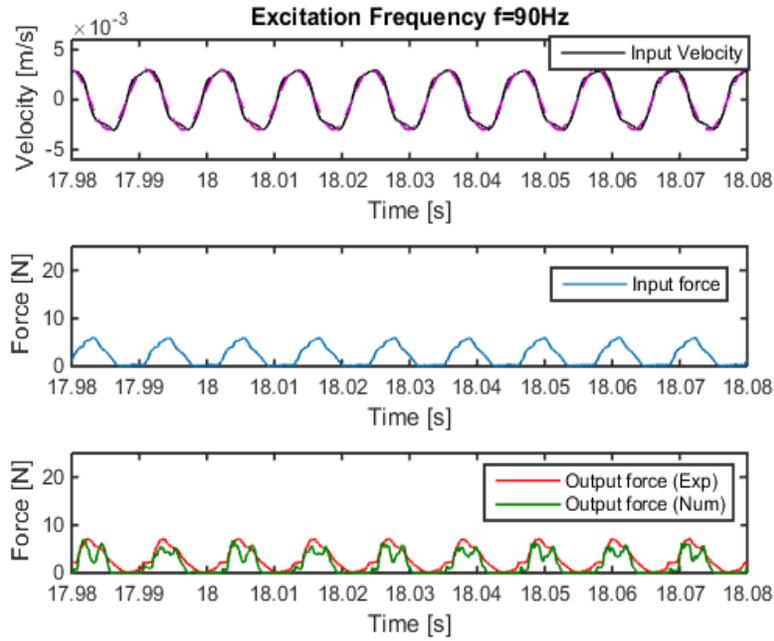

(a)

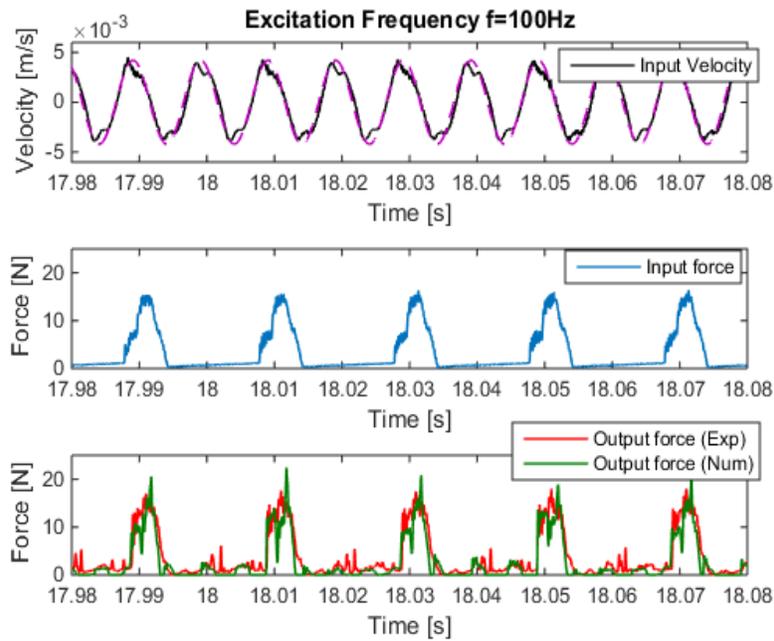

(b)

Figure 10. Comparisons between experimental measurements (red) and numerical results (green) from the model (3): Experimental velocity time series of the armature of the shaker (*top*), experimental input force at the left end of the first granule (*middle*) and transmitted force at the right end (*bottom*) for (a) 1:1 anti-resonance at 90 Hz, and (b) 1:2 resonance at 100 Hz.

In Figure 9 we compare the experimental transmitted forces directly recorded by the force transducer with direct numerical simulations of the granular network (3) with a fixed boundary condition at the end; in the computational model the output force is calaculated as $F_{out} = 2E\sqrt{R}(u_{11})_+^{3/2}/3(1-v^2)$. In Figure 9a the maximum of the numerical output force at the right end of the chain is superimposed on the experimental results in the frequency range $20Hz \leq f \leq 160Hz$. Satisfactory agreement is inferred, with the computational predictions fully capturing the two resonances and the anti-resonance in-between of the experimental measurements. Moreover, comparisons of the experimental and computational time series of the transmitted force at 60 Hz (1:1 resonance), 90 Hz (1:1 anti-resonance) and 100 Hz (1:2 resonance) are depicted in Figs. 9b, 10a and 10b, respectively. Examining the force pulse trains depicted in these plots, we clearly deduce that all numerical simulations reproduce accurately the experimental measurements, following the same trends and even reproducing the ringing effects due to the flexural supports at the two resonances (cf. Figs. 9b and 10b). It follows that the augmented mathematical model (3) is capable of accurately reproducing the experimental measurements, fully validating the experimentally detected strongly nonlinear resonance and anti-resonance steady state responses.

The deviations between the experimental and computational results may be attributed to the approximation of dissipative effects by linear viscous damping used in the numerical model (3) which is incapable of fully modeling nonlinear dissipative effects such as friction and plasticity. Despite numerous efforts [1, 21, 28], a universal model capturing quantitatively the phenomenology of dissipative losses is still not available. In addition, possible inherent bead misalignments in the experimental fixture can affect the force pulses transmitted to the right end of the granular chain, as well as the theoretical modeling of the right boundary condition as fixed, i.e., of infinite stiffness. It is clear that in the practical realization of the granular chain, the right boundary is a force transducer which has a finite stiffness. Nevertheless, the computational results are in satisfactory agreement with the experimental measurements, both in the frequency and the time domains, fully recovering the resonances and anti-resonances in the stationary-state responses.

## 4. Data assimilation and stability analysis

In this section we provide additional justification for the inclusion of the weak foundation from the flexures attached to each bead by conducting a numerical stability analysis along with validation of experimental results using a numerical model similar to (3). As mentioned before, the system in (3) is provided with an experimental time series of input force and subsequently integrated in time to obtain the numerical results. In effect, these results are nonlinear observers of the system state, obtained from assimilation of the experimental force data. In order to study the numerical stability of the system response we construct the following extended system based on the dimensional equations in (3).

$$\ddot{x}_1 = \left(\sin(\beta\tau) - x_1\right)_+^{3/2} - \left(x_1 - x_2\right)_+^{3/2}$$
$$+ \lambda\left\{\left(\beta\cos(\beta\tau) - \dot{x}_1\right)H\left(\sin(\beta\tau) - x_1\right) - \left(\dot{x}_1 - \dot{x}_2\right)H\left(x_1 - x_2\right)\right\} - kx_1$$
$$\ddot{x}_i = \left(x_{i-1} - x_i\right)_+^{3/2} - \left(x_i - x_{i+1}\right)_+^{3/2}$$
$$+ \lambda\left\{\left(\dot{x}_{i-1} - \dot{x}_i\right)H\left(x_{i-1} - x_i\right) - \left(\dot{x}_i - \dot{x}_{i+1}\right)H\left(x_i - x_{i+1}\right)\right\} - kx_i, \quad i = 2,\ldots,N-1$$
$$\ddot{x}_N = \left(x_{N-1} - x_N\right)_+^{3/2} - \sqrt{2}\left(x_N\right)_+^{3/2} + \lambda\left\{\left(\dot{x}_{N-1} - \dot{x}_N\right)H\left(x_{N-1} - x_N\right) - \left(\dot{x}_N\right)H\left(x_N\right)\right\} - kx_N; \quad (4)$$
$$\ddot{y}_1 = F_{in} - \left(y_1 - y_2\right)_+^{3/2} + \lambda\left\{-\left(\dot{y}_1 - \dot{y}_2\right)H\left(y_1 - y_2\right)\right\} - ky_1$$
$$\ddot{y}_i = \left(y_{i-1} - y_i\right)_+^{3/2} - \left(y_i - y_{i+1}\right)_+^{3/2}$$
$$+ \lambda\left[\left(\dot{y}_{i-1} - \dot{y}_i\right)H\left(y_{i-1} - y_i\right) - \left(\dot{y}_i - \dot{y}_{i+1}\right)H\left(y_i - y_{i+1}\right)\right] - ky_i, \quad i = 2,\ldots,N-1$$
$$\ddot{y}_N = \left(y_{N-1} - y_N\right)_+^{3/2} - \sqrt{2}\left(y_N\right)_+^{3/2} + \lambda\left[\left(\dot{y}_{N-1} - \dot{y}_N\right)H\left(y_{N-1} - y_N\right) - \left(\dot{y}_N\right)H\left(y_N\right)\right] - ky_N,$$
$$F_{in} = \left(\sin(\beta\tau) - x_1\right)_+^{3/2} + \lambda\left(\beta\cos(\beta\tau) - \dot{x}_1\right)H\left(\sin(\beta\tau) - x_1\right)$$

where $k = 0.001$ is the dimensionless foundation stiffness. The extended model (4) is in non-dimensional form incorporating the same normalizations used to derive equation (2). These equations correspond to two granular chains where the first one is independent from the second and periodically forced on its left end, while the second granular chain corresponds to the model (3), yet its forcing ($F_{in}$) is prescribed from a different source, namely, the input signal from the first granular chain. This is again a data assimilation problem: one uses the second chain equations, forced by the first chain, to create nonlinear observers of the first chain states based on measurements of force on the leftmost bead of this chain. The solutions of the above equations are periodic orbits in phase space, which correspond to fixed points of the stroboscopic (Poincaré) map with period $1/f$. Applying Newton's method to solve the nonlinear boundary value problem we are able to obtain the Floquet multipliers that show the stability of the overall coupled system (the first chain forcing the second chain).

Two cases are considered, one where the weak foundation has been omitted by using $k = 0$ and one that implements the value $k = 0.001$. The Floquet multipliers of the combined system of chains, along with the displacement of the first bead in each granular chain are depicted in Figure 11. The top row of panels corresponds to the case of $k = 0$, and shows that there is one *unstable* multiplier with modulus greater than one; this means that while the periodic orbit of the first granular chain, forced by a periodically moving $0^{th}$ bead is stable, the periodic orbit of the second chain, taking as input the force on its first bead the corresponding force from the first chain, will eventually "break down" (that is, its linear instability will become manifest by the solution moving away from it). The data assimilation formulation is then *unstable*, the models *do not synchronize* [12, 17, 27] and we do not obtain good observations of the experimental state from the model and the input forces. This shows that a granular chain without grounding stiffnesses (modeling the effects of the supporting flexures) that is forced externally by using an experimental input force is inherently unstable and will not lead to time periodic solutions under these conditions. Including the weak foundation from the flexures makes the system numerically stable as becomes evident in the bottom row of panels of Fig. 11. Every multiplier is now inside the unit circle, the periodic solution for both chains is stable, and the displacements of the first bead of each chain are identical since the input signal to both granular chains is the same.

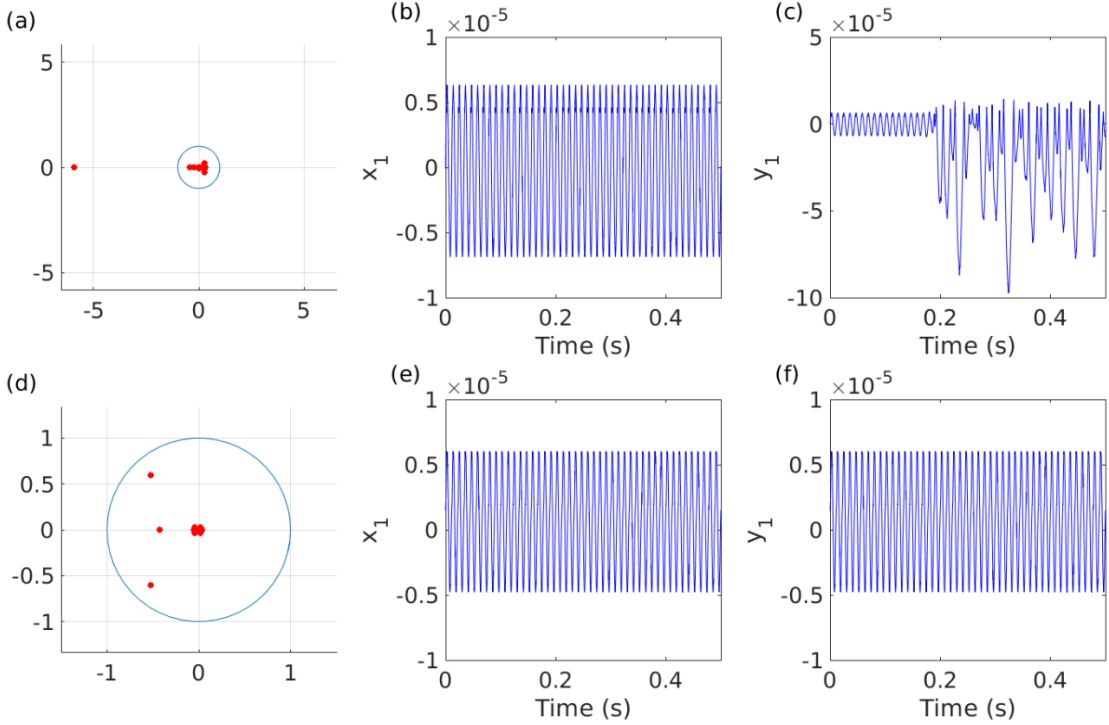

Figure 11. Stability analysis of the extended model (4): (a) One unstable Floquet multiplier when the weak foundation is omitted, (b) the periodic solution of the independent granular chain is stable, (c) the periodic solution of the second granular chain is unstable and eventually "breaks down"; (d) Floquet multiplier spectrum for $k = 0.001$, the model is now stable, (e,f) the periodic solutions of both granular chains are stable and identical since the input signal is the same.

We now proceed with simulating the experimental results using a model similar to the one in (3). We substitute, however, the $F_{in}$ term with the force between the first bead and the shaker, which we model as a $0^{th}$ bead with prescribed harmonic displacement, $u_0 = A_0 \sin(2\pi f t)$, or in dimensionless form, $y_0 = \sin(\beta\tau)$. The equation of the $1^{st}$ bead of the second chain becomes:

$$\ddot{y}_1 = (\sin(\beta\tau) - y_1)_+^{3/2} - (y_1 - y_2)_+^{3/2}$$
$$+ \lambda\{(\beta\cos(\beta\tau) - \dot{y}_1)H(\sin(\beta\tau) - y_1) - (\dot{y}_1 - \dot{y}_2)H(y_1 - y_2)\} - ky_1 \quad (5)$$

The physical parameters of the system remain the same as in (3) and we adjust the amplitude of the external forcing and compare to the experimental results. Figure 12 shows the results at the frequencies of 60Hz (1:1 resonance), 90 Hz (1:1 anti-resonance) and 100 Hz (1:2 resonance). The top row of the panels depicts the velocity time series of the shaker both in the model and in the experiment, the middle row shows the input force time series on the left end of the chain, and the bottom row shows the transmitted force time series at the right end of the chain.

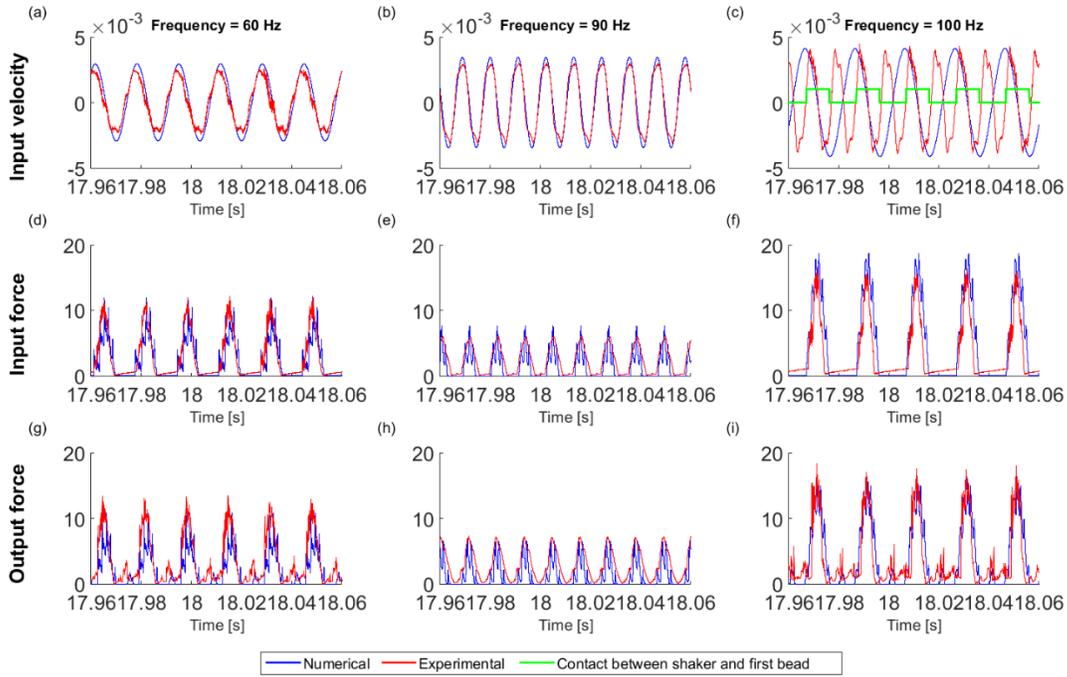

Figure 12. Results of numerical modeling using a 0th bead with prescribed harmonic displacement. (a),(b),(c) The velocity of the shaker is replicated accurately at 60 and 90 Hz as well as at 100 Hz when the shaker is actually in contact with the granular chain. (d),(e),(f) The input force at the left end of the chain. (g),(h),(i) The transmitted force at the right end of the chain. Both the strong pulses and the weak pulses are present in the numerical results as well.

It is clear that the numerical results at 60 Hz and 90 Hz match the experimental ones quite accurately, including the small pulses ("ringing" due to the dynamics of the supporting flexures) of the output force as seen in Figure 12g. However, in order to replicate the experiment at 100 Hz, we had to use a modified frequency $f$ = 50 Hz in our model. The reason for this discrepancy is that while the shaker in the experiment is driven at 100 Hz, a hit occurs every *second* period of its motion making it a 1:2 resonance with a period-2 harmonic solution. The green line in Figure 12c is a Heaviside function, showing when the shaker is in contact with the first bead in the numerical model. As we can see, when force is being applied to the chain, the experimental velocity of the shaker is replicated accurately (in the "eye" norm) by our model. Additionally, the small pulses of the output force in Figure 12i were replicated as well. The actual solution of the model at 100 Hz was a period-1 solution but this difference can be attributed to the coexistence of two stable solutions at the same frequency and the convergence of the experiment to the period-2 solution. One caveat of this modeling approach is that different values of $A_0$ had to be used at each frequency. This is corroborated by the fact that the experimentally computed amplitude also lies within a range of ~6 − 8 μm. The range of amplitudes used in the numerical model, ~6 − 13 μm, while larger is still in reasonably good agreement with the experimental results.

## 5. Conclusions

In this work, we have studied the stationary-state dynamics of a one-dimensional finite homogeneous granular chain, with no prior compression and under time-periodic excitation. We experimentally confirmed two types of resonance motions involving harmonic or subharmonic traveling pulses in the granular chain, as well as a state of anti-resonance. These results, which correspond to local maxima and minima, respectively, of the maximum transmitted force at the right end of the chain, validate previously reported theoretical studies. In particular, in agreement with previous theoretical predictions we experimentally verified the existence of two strongly nonlinear resonance peaks and one anti-resonance valley between them within the frequency range of $20 - 160\ Hz$. Hence, we experimentally proved that a strongly nonlinear medium (with complete absence of linear acoustics and with no linear resonance spectrum) can still support a nonlinear resonance spectrum that is tunable with energy. To the authors' knowledge this is the first such experimental result reported in the literature. Furthermore, we revisited the mathematical model of the experimental fixture and were able to confirm the existence of aforementioned interesting dynamic responses by direct numerical simulations.

From the results obtained with the experimental granular chains corresponding to three frequencies ($60\ Hz$, $90\ Hz$ and $100\ Hz$) we ascertained that the theoretically predicted resonances and anti-resonances can be realized experimentally. The findings represented in this work have multiple potential applications in the design of acoustic metamaterials which are strongly nonlinear. These media can be used as energy absorbers when they are excited at anti-resonance frequencies, but, on the contrary, can intensify energy transmission at resonance frequencies. We emphasize that due to the strong nonlinearity of the granular dynamics, the detected resonance spectra are dependent on the intensities (magnitudes) of the applied excitations, so the stationary-state dynamics are passively tunable with energy. Moreover, such designs can be extended to higher dimensions, exploring the concepts of resonance and anti-resonance in these setups as well, e.g., in hexagonal or square lattices. In general, the results of this work contribute to the design of practical nonlinear acoustic metamaterials with properties adaptive to different types of external excitations.

## 6. Acknowledgments

YZ, DMM, and AFV gratefully acknowledge the support of US ARO through MURI grant W911NF-09-1-0436. Dr. David Stepp was the grant monitor. DP, PGK and IGK would like to acknowledge the support of US AFOSR through grant FA9550-12-1-0332. PGK also acknowledges the support of ARO (W911NF-15-1-0604). IGK also acknowledges partial support by the US National Science Foundation (grant no. ECCS 1406224).

# References

1. Carretero-Gonzalez R, Khatri D, Porter MA, Kevrekidis PG, Daraio C (2009) Dissipative solitary waves in granular crystals. Phys. Rev. Lett., 102, 024102

2. Charalampidis EG, Li F, Chong C, Yang J, Kevrekidis PG (2015) Time-periodic solutions of driven damped trimer granular crystals. Mathematical Problems in Engineering, 830978

3. Coste C, Falcon E, Fauve S (1997) Solitary waves in a chain of beads under Hertz contact. Phys. Rev. E, 56 (5), 6104–6117

4. Chong C, Li F, Yang J, Williams MO, Kevrekidis IG, Kevrekidis PG, Daraio C (2014) Damped-driven granular chains: An ideal playground for dark breathers and multibreathers. Phys. Rev. E, 89(3), 032924

5. Daraio C, Nesterenko VF, Herbold B, Jin S (2006) Tunability of Solitary Wave Properties in One-dimensional Strongly Nonlinear Phononic Crystals. Phys. Rev. E, 73(2), 026610.

6. Hasan M, Cho S, Remick K, Vakakis AF, McFarland DM, Kriven W (2015) Experimental Study of Nonlinear Acoustic Bands and Propagating Breathers in Ordered Granular Media Embedded in Matrix. Gran. Matter, 17, 49–72

7. Herbold EB, Nesterenko VF (2007) Shock wave structure in a strongly nonlinear lattice with viscous dissipation. Phys. Rev. E, 75(2), 021304.

8. Hoogeboom C, Man Y, Boechler N, Theocharis G, Kevrekidis PG, Kevrekidis IG, Daraio C (2013) Hysteresis loops and multi-stability: From periodic orbits to chaotic dynamics (and back) in diatomic granular crystals. EPL, 101, 44003

9. James G (2011) Nonlinear Waves in Newton's Cradle and the Discrete p-Schrödinger Equation. Math. Mod. Meth. Appl. Sci., 21, 2335-2377

10. Jayaprakash KT, Starosvetsky Y, Vakakis AF, Peeters M, Kerschen G (2011) Nonlinear Normal Modes and Band Zones in Granular Chains with No Precompression. Nonl. Dyn., 63(3), 359-385

11. Kevrekidis PG (2011) Non-linear waves in lattices: past, present, future. IMA J. Appl. Math., 76, 389-423

12. Law KJH, Stuart AM, Zygalakis KC (2015) Data Assimilation: A Mathematical Introduction, Springer

13. Lazaridi AN, Nesterenko VF (1985) Observation of a New Type of Solitary Waves in a One-Dimensional Granular Medium. J. Appl. Mech. Tech. Physics, 26(3), 405-408

14. Lydon J, Jayaprakash KR, Ngo D, Starosvetsky Y, Vakakis AF, Daraio C (2013) Frequency Bands of Strongly Nonlinear Finite Homogeneous Granular Chains. Phys. Rev. E, 88, 012206

15. Nesterenko VF (2001), Dynamics of Heterogeneous Materials, Springer Verlag, New York

16. Nesterenko VF (1983) Propagation of Nonlinear Compression Pulses in Granular Media. J. Appl. Mech. Techn. Physics, 24(5), 733-743

17. Pecora LM, Carroll TL (1990) Synchronization in chaotic systems. Phys. Rev. Lett., 64(8), 821

18. Porter MA, Kevrekidis PG, Daraio C (2015) Granular Crystals: Nonlinear dynamics meets materials engineering. Physics Today, 68, 44-50

19. Potekin R, Jayaprakash KR, McFarland DM, Remick K, Bergman LA, and Vakakis AF (2012) Experimental study of strongly nonlinear resonances and anti-resonances in granular dimer chains. Exp. Mech., 53(5), 861–870.

20. Pozharskiy D, Zhang Y, Williams M, McFarland DM, Kevrekidis PG, Vakakis AF, Kevrekidis IG (2015) Nonlinear Resonances and Antiresonances of a Forced Sonic Vacuum. Phys. Rev. E, 92(6), 063203; and arXiv:1507.01025

21. Rosas A, Romero AH, Nesterenko VF, Lindenberg K (2007) Observation of two-wave structure in strongly nonlinear dissipative granular chains. Phys. Rev. Lett., 98, 164301

22. Sen S, Manciu M, Wright JD (1998) Soliton-like pulses in perturbed and driven Hertzian chains and their possible applications in detecting buried impurities. Phys. Rev. E, 57 (2), 2386–2397

23. Sen S, Hong J, Bang J, Avalos E, Doney R (2008) Solitary waves in the granular chain. Phys. Rep., 462 (2), 21–66

24. Sinkovits RS, Sen S (1995) Nonlinear dynamics in granular columns. Phys. Rev. Lett., 74 (14), 2686–2689

25. Spadoni A, Daraio C (2010) Generation and control of sound bullets with a nonlinear acoustic lens. Proc. Natl. Acad. Sci., 107 (16), 7230–7234

26. Starosvetsky Y, Vakakis AF (2010) Traveling Waves and Localized Modes in One-Dimensional Homogeneous Granular Chains With no Pre-Compression. Phys. Rev. E, 82(2), 026603

27. Strogatz SH (2001) Nonlinear dynamics and chaos with applications to physics, biology, chemistry, and engineering, Westview press

28. Vergara L (2010) Model for dissipative highly nonlinear waves in dry granular crystals. Phys. Rev. Lett., 104, 244302

29. Yang J, Sutton M (2015) Nonlinear wave propagation in hexagonally packed granular channel under rotational dynamics. Int. J. Solids Str., doi: 10.1016/j.ijsolstr.2015.07.017

30. Zhang Y, Moore KJ, McFarland DM, Vakakis AF (2015) Targeted Energy Transfers and Passive Acoustic Wave Redirection in a Two-dimensional Granular Network Under Periodic Excitation. J. Appl. Phys., 118, 234901